\documentclass[tikz, sepfignums]{nmd/article}
\ifnmd
\fi
\usepackage{enumitem}
\usepackage{colortbl}
\usepackage{placeins}
\usetikzlibrary{matrix, shapes.geometric}

\title{A census of exceptional \\ Dehn fillings}

\author{Nathan M. Dunfield}
\givenname{Nathan}
\surname{Dunfield}
\address{ Dept.~of Math., MC-382 \\
          University of Illinois \\
          1409 W. Green St. \\
          Urbana, IL 61801 \\ 
          USA
}
\email{nathan@dunfield.info}
\urladdr{http://dunfield.info}

\arxivreference{}
\arxivpassword{}

\newcommand{\slopes}{\mathit{Sl}}

\newcommand{\connsum}{\mathop{\#}}

\pgfkeys{/matplotlibfigure,
  default/.style = {
    width = 0.45\textwidth,
    font = \fontsize{9}{10.7}\selectfont}
}

\begin{document}

\begin{abstract} 
  This paper describes the complete list of all 205,822 exceptional
  Dehn fillings on the 1-cusped hyperbolic 3-manifolds that have ideal
  triangulations with at most 9 ideal tetrahedra.  The data is
  consistent with the standard conjectures about Dehn filling and
  suggests some new ones.
\end{abstract}
\maketitle

\section{Introduction}

\subsection{Dehn filling} Suppose $M$ is a compact orientable
\3-manifold with $\partial M$ a torus.  A \emph{slope} on $\partial M$
is an unoriented isotopy class of simple closed curve, or equivalently
a primitive element of $H_1(\partial M; \Z)$ modulo sign.  The set of
all slopes will be denoted $\slopes(M)$, which can be viewed as the
rational points in the projective line
$P^1\big(H_1(\partial M; \R)\big) \cong \RP^1$. The Dehn fillings of
$M$ are parameterized by $\alpha \in \slopes(M)$, with $M(\alpha)$
being the Dehn filling where $\alpha$ bounds a disk in the attached
solid torus.  When the interior of $M$ admits a hyperbolic metric of
finite volume, it is called a \emph{\1-cusped hyperbolic
  \3-manifold}. For such hyperbolic $M$, Thurston showed that all but
finitely many $M(\alpha)$ are also hyperbolic \cite{ThurstonsNotes}.
The nonhyperbolic Dehn fillings are called \emph{exceptional}, and the
corresponding slopes the \emph{exceptional slopes}.  Understanding the
possible exceptional fillings has been a major topic in the study of
3-manifolds over the past 40 years; see the surveys \cite{Gordon1998,
  Gordon1999, Gordon2003, Gordon2012} for further background.

\begin{figure}
  \newcommand{\countbytettable}{%
    \footnotesize
    \begin{tabular}{rr}
      \toprule
    {tets} & \#manifolds \\
    \midrule
      2    &         2 \\
      3    &         9 \\
      4    &        52 \\
      5    &       223 \\
      6    &       913 \\
      7    &   3{,}388 \\
      8    &  12{,}241 \\
      9    &  42{,}279 \\
      \midrule
      total & 59{,}107 \\
      \toprule
    \end{tabular}
  }
  \begin{center}
    {%
 \pgfkeys{/matplotlibfigure, default, }%
 \begin{tikzoverlay}[width=\matplotlibfigurewidth]{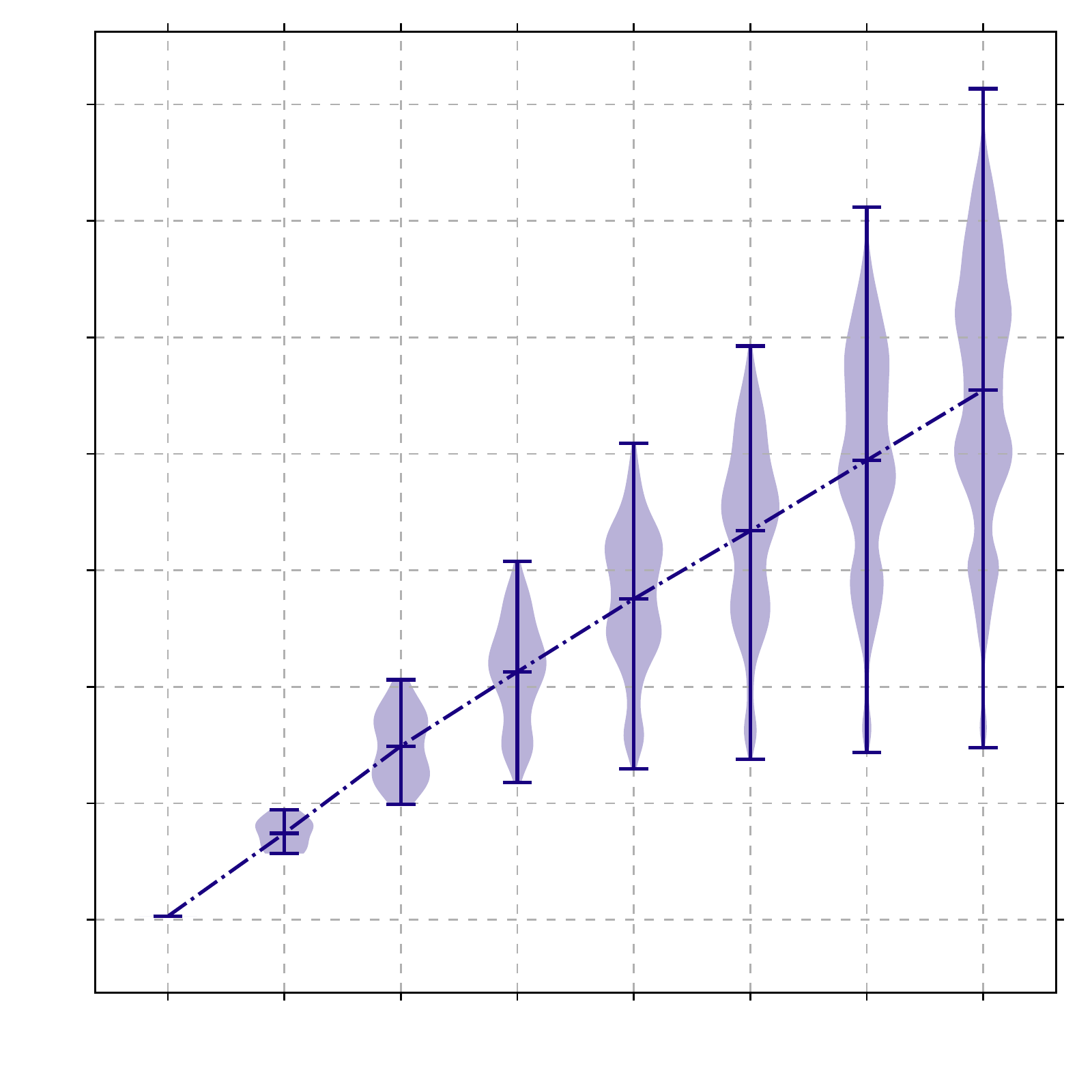}[\matplotlibfigurefont]
  \draw (52.707248, 3.326270) node[below=0.4em] {\footnotesize tetrahedra};
  \draw (15.372442, 7.586806) node[below] {$2$};
  \draw (26.039529, 7.586806) node[below] {$3$};
  \draw (36.706617, 7.586806) node[below] {$4$};
  \draw (47.373704, 7.586806) node[below] {$5$};
  \draw (58.040792, 7.586806) node[below] {$6$};
  \draw (68.707880, 7.586806) node[below] {$7$};
  \draw (79.374967, 7.586806) node[below] {$8$};
  \draw (90.042055, 7.586806) node[below] {$9$};
  \draw (-2, 53.107639) node[rotate=90.0] {\footnotesize volume};
  \draw (7.186415, 15.772832) node[left] {$2$};
  \draw (7.186415, 26.439920) node[left] {$3$};
  \draw (7.186415, 37.107008) node[left] {$4$};
  \draw (7.186415, 47.774095) node[left] {$5$};
  \draw (7.186415, 58.441183) node[left] {$6$};
  \draw (7.186415, 69.108270) node[left] {$7$};
  \draw (7.186415, 79.775358) node[left] {$8$};
  \draw (7.186415, 90.442445) node[left] {$9$};
  \begin{scope}[shift={(-5.96173322, -5.56134259)},
                xscale=10.66708754, yscale=10.66708754]
  \end{scope}
\end{tikzoverlay}
 \hspace{0.5em} %
 \pgfkeys{/matplotlibfigure, default, }%
 \begin{tikzoverlay}[width=\matplotlibfigurewidth]{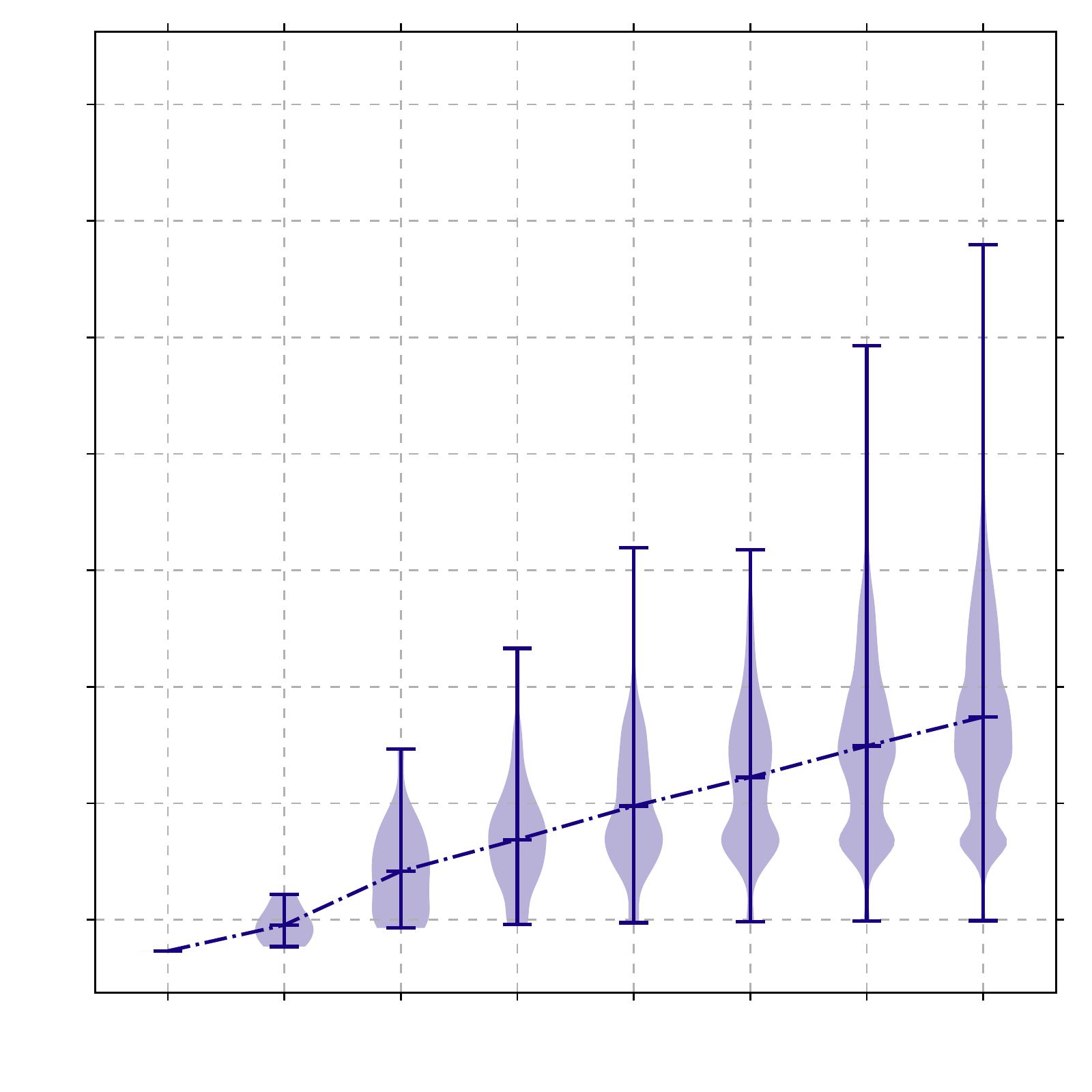}[\matplotlibfigurefont]
  \draw (52.707248, 3.326270) node[below=0.4em] {\footnotesize tetrahedra};
  \draw (15.372442, 7.586806) node[below] {$2$};
  \draw (26.039529, 7.586806) node[below] {$3$};
  \draw (36.706617, 7.586806) node[below] {$4$};
  \draw (47.373704, 7.586806) node[below] {$5$};
  \draw (58.040792, 7.586806) node[below] {$6$};
  \draw (68.707880, 7.586806) node[below] {$7$};
  \draw (79.374967, 7.586806) node[below] {$8$};
  \draw (90.042055, 7.586806) node[below] {$9$};
  \draw (-2, 53.107639) node[rotate=90.0] {\footnotesize cusp volume};
  \draw (7.186415, 15.772832) node[left] {$2$};
  \draw (7.186415, 26.439920) node[left] {$3$};
  \draw (7.186415, 37.107008) node[left] {$4$};
  \draw (7.186415, 47.774095) node[left] {$5$};
  \draw (7.186415, 58.441183) node[left] {$6$};
  \draw (7.186415, 69.108270) node[left] {$7$};
  \draw (7.186415, 79.775358) node[left] {$8$};
  \draw (7.186415, 90.442445) node[left] {$9$};
  \begin{scope}[shift={(-5.96173322, -5.56134259)},
                xscale=10.66708754, yscale=10.66708754]
  \end{scope}
\end{tikzoverlay}
}
    
    \hspace{3em} \begin{minipage}{0.25\textwidth}
      \countbytettable
      \vspace{3em}
    \end{minipage}
    \hspace{1.6em}
    \begin{minipage}{0.55\textwidth}
 \pgfkeys{/matplotlibfigure, default, width=\textwidth}%
 \begin{tikzoverlay}[width=\matplotlibfigurewidth]{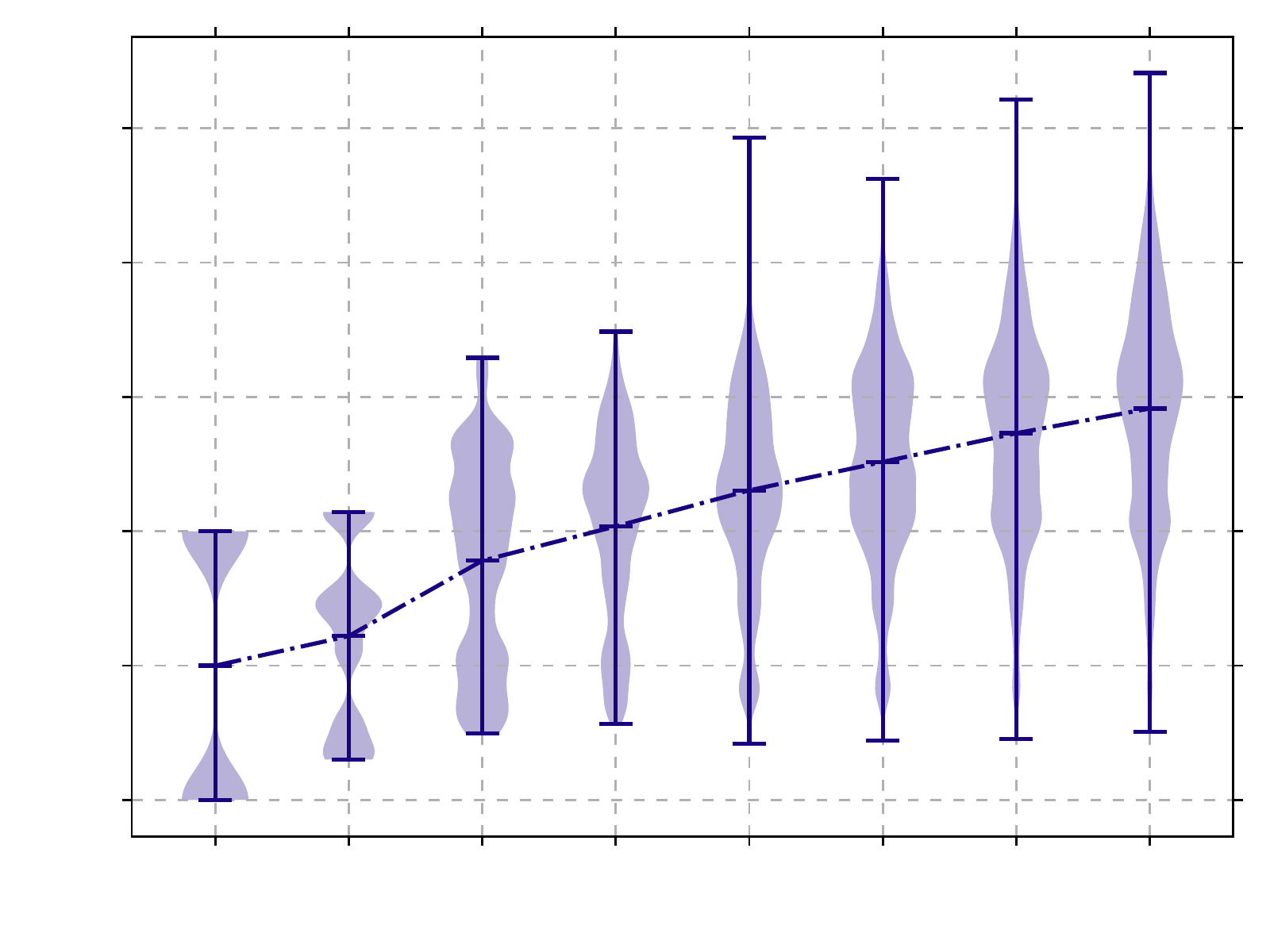}[\matplotlibfigurefont]
  \draw (53.742405, 4.264874) node[below=0.4em] {\footnotesize tetrahedra};
  \draw (16.946187, 7.586806) node[below] {$2$};
  \draw (27.459392, 7.586806) node[below] {$3$};
  \draw (37.972597, 7.586806) node[below] {$4$};
  \draw (48.485802, 7.586806) node[below] {$5$};
  \draw (58.999007, 7.586806) node[below] {$6$};
  \draw (69.512212, 7.586806) node[below] {$7$};
  \draw (80.025417, 7.586806) node[below] {$8$};
  \draw (90.538622, 7.586806) node[below] {$9$};
  \draw (-2, 37.9) node[rotate=90.0] {\footnotesize min slope length};
  \draw (8.856337, 11.969697) node[left] {$1.0$};
  \draw (8.856337, 22.556563) node[left] {$1.5$};
  \draw (8.856337, 33.143428) node[left] {$2.0$};
  \draw (8.856337, 43.730294) node[left] {$2.5$};
  \draw (8.856337, 54.317160) node[left] {$3.0$};
  \draw (8.856337, 64.904026) node[left] {$3.5$};
  \begin{scope}[shift={(-4.08022280, -9.20403452)},
                xscale=10.51320497, yscale=21.17373149]
  \end{scope}
\end{tikzoverlay}

    \end{minipage}
  \end{center}

  \vspace{-1em}
  
  \caption{Some basic statistics on the manifolds in $\cC_9$.  The
    table in the lower left shows the number of manifolds in
    $\cC_t \setminus \cC_{t-1}$, that is, those whose minimal ideal
    triangulations have exactly $t$ tetrahedra. In the top left, the
    hyperbolic volume is shown via a violin plot.  Here, each
    ``violin'' shows the distribution of volume for
    $\cC_t \setminus \cC_{t-1}$, where the top and bottom horizontal
    bars are the min and max, the middle horizontal bar is the mean,
    and the violin ``body'' is a smoothed histogram of the volumes.
    The plot in the upper right shows the cusp volume, i.e.~the
    volume of a maximal cusp neighborhood bounded by an embedded
    horotorus.  The plot in the lower right shows the minimal slope
    length, that is, the shortest essential curve in the maximal
    horotorus. See Section~\ref{sec:slopelength} for the relevance of
    this cusp data.}
  \label{fig:cusped}
\end{figure}

\begin{figure}
  \renewcommand\figurename{Table}
  \begin{center}
  \footnotesize
  \begin{tabular}{p{0.1cm}lr@{\hskip 1.5cm}lr}
    \toprule
    \multicolumn{3}{c}{$b_1\big(M(\alpha)\big)$  = 0} & \multicolumn{2}{c}{$b_1\big(M(\alpha)\big)$  > 0}\\
    \midrule
    \multicolumn{2}{l}{atoroidal}                        \\ 
    & $S^3$           &  1{,}267                         \\ 
    & lens space      & 44{,}487                       \\ 
    & finite $\pi_1$  & 13{,}446 & $S^2 \times S^1$ & 242 \\
    & Seifert fibered & 71{,}111 & Seifert fibered & 118 \\ 
    & connected sum   &  4{,}296 & connected sum   & 169 \\
    \\ 
    \multicolumn{2}{l}{toroidal} \\ 
    & Seifert fibered &  1{,}730 & Seifert fibered   &     159 \\ 
    & graph manifold  & 63{,}325 & graph manifold    & 3{,}043 \\
    & hyperbolic piece&  2{,}136 & hyperbolic piece  &      74 \\
    &                 &          & Sol torus bundle &     219 \\
    \midrule
    \multicolumn{2}{l}{totals} & 201{,}798 & & 4{,}024\\
    \bottomrule
  \end{tabular}
\end{center}

\caption{Summary of the topological types of the 205{,}822 exceptional
  fillings from Theorem~\ref{thm:except}, broken down by betti number
  and whether the filling is toroidal. See Section~\ref{sec:conven}
  for precise definitions. Here, each filling is listed only once in
  the most restricted category possible from
  Figure~\ref{fig:excepvenn}, and there are no fillings that are both
  toroidal and connected sums.}
\label{table:exceptsum}

  \vspace*{\floatsep}

  \renewcommand\figurename{Figure}  
\newsavebox{\eMtablebox}
\savebox{\eMtablebox}{%
  \footnotesize
  \begin{tabular}{rr}
    \toprule
    $e(M)$ & \#manifolds\\
    \midrule
    0   &         3 \\
    1   &   1{,}824 \\
    2   &   9{,}730 \\
    3   &  18{,}003 \\
    4   &  19{,}009 \\
    5   &   8{,}751 \\
    6   &   1{,}776 \\
    7   &         8 \\
    8   &         2 \\
    9   &         0 \\
    10  &         1 \\
    \bottomrule
  \end{tabular}
}

\hspace{1.0em} 
{
\newlength\eMtableheight
\setlength\eMtableheight{\dimexpr \ht\eMtablebox+\dp\eMtablebox}
\pgfkeys{/matplotlibfigure, default}
\begin{tikzoverlay}[height=1.188\eMtableheight]{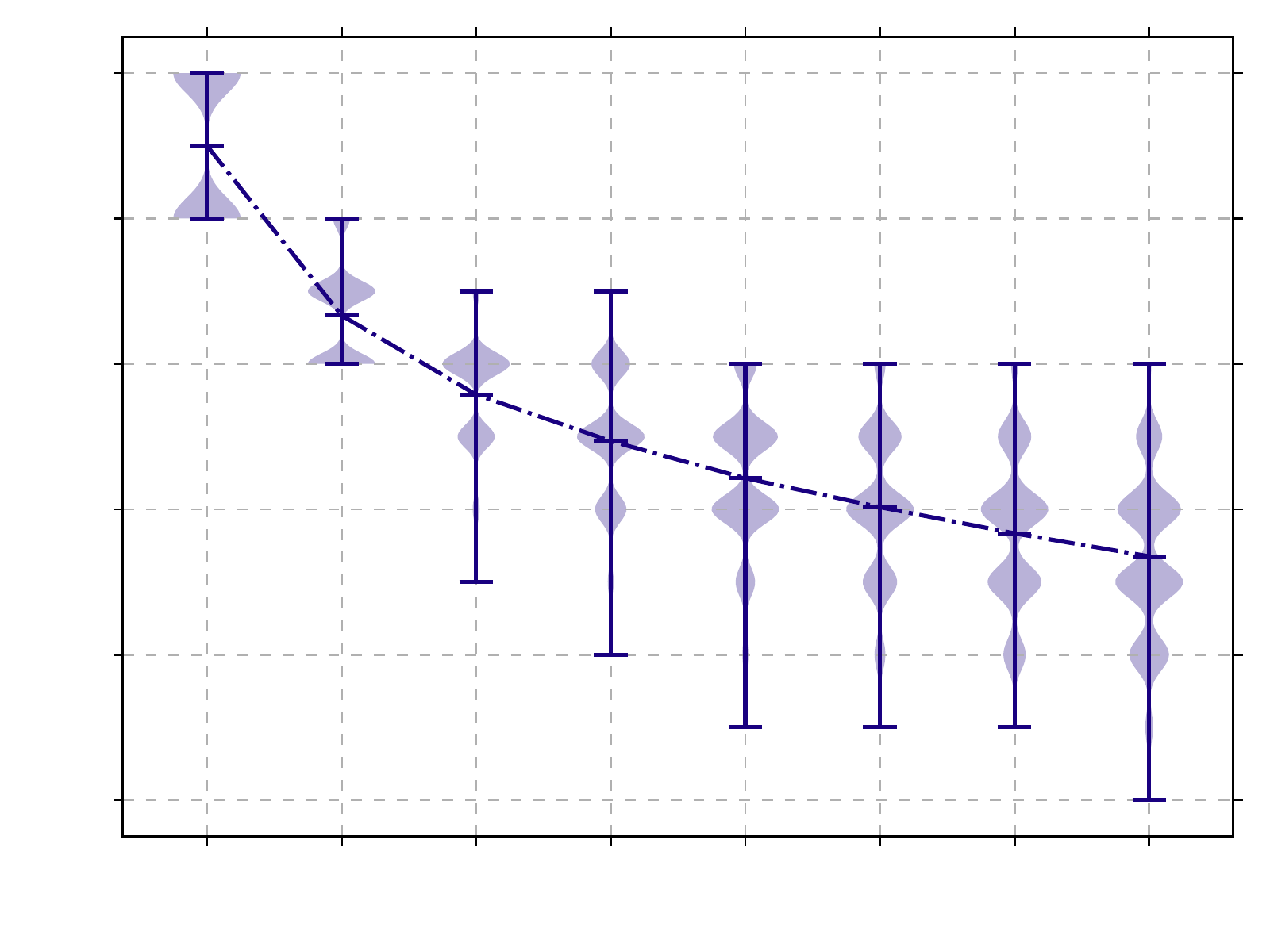}[font=\matplotlibfigurefont]
  \node[anchor=south east, inner sep=0] at (-16, 9.0) {\usebox{\eMtablebox}};
  %
  %
  \draw (53.390842, 3.326270) node[below=0.4em] {\footnotesize tetrahedra};
  \draw (16.296329, 7.586806) node[below] {$2$};
  \draw (26.894761, 7.586806) node[below] {$3$};
  \draw (37.493194, 7.586806) node[below] {$4$};
  \draw (48.091626, 7.586806) node[below] {$5$};
  \draw (58.690058, 7.586806) node[below] {$6$};
  \draw (69.288490, 7.586806) node[below] {$7$};
  \draw (79.886923, 7.586806) node[below] {$8$};
  \draw (90.485355, 7.586806) node[below] {$9$};
  \coordinate (eM) at (3.258008, 40.607639);
  \node[left] at (eM) {\footnotesize $e(M)$};
  \draw (8.153212, 11.969697) node[left] {$0$};
  \draw (8.153212, 23.424874) node[left] {$2$};
  \draw (8.153212, 34.880051) node[left] {$4$};
  \draw (8.153212, 46.335227) node[left] {$6$};
  \draw (8.153212, 57.790404) node[left] {$8$};
  \draw (8.153212, 69.245581) node[left] {$10$};
  \begin{scope}[shift={(-4.90053530, 11.96969697)},
                xscale=10.59843224, yscale=5.72758838]
  \end{scope}
\end{tikzoverlay}
}

\caption{This figure describes how the exceptional fillings of
  Theorem~\ref{thm:except} are distributed over the manifolds in
  $\cC_9$.  The number of exceptional fillings on $M$ is denoted
  $e(M)$, and the table at left shows the number of manifolds in
  $\cC_9$ having each possible value of $e(M)$, which is at most
  $10$ by \cite{LackenbyMeyerhoff2013}. At right is a violin plot of
  $e(M)$ as a function of the number of tetrahedra of $M$; see
  Figure~\ref{fig:cusped} for more on violin
  plots.}
\label{fig:numexcep}

\end{figure}

This paper gives a census of all exceptional Dehn fillings on a
certain collection of \1-cusped hyperbolic \3-manifolds.
Specifically, let $\cC_t$ be the set of all orientable \1-cusped
hyperbolic \3-manifolds that have ideal triangulations with at most
$t$ ideal tetrahedra.  For $t \leq 9$, the set $\cC_t$ has been
enumerated by \cite{HildebrandWeeks1989, CallahanHildebrandWeeks1999,
  Thistlethwaite2010, Burton2014} and is included with SnapPy
\cite{SnapPy}, whose nomenclature for these manifolds (e.g.~$m004$,
$s011$, $v1002$, $t12345$, and $o9_{60000}$) I will use freely
throughout. Each manifold $M$ in $\cC_t$ has a preferred basis for
$H_1(\partial M; \Z)$, and so I will denote slopes in $\slopes(M)$ by
elements in $\Z^2$.  See Figure~\ref{fig:cusped} for some basic
statistics on the $59{,}107$ manifolds in $\cC_9$.  The main result of
this paper is:
\begin{theorem}\label{thm:except}
  There are precisely 205{,}822 exceptional Dehn fillings on the
  manifolds in $\cC_9 $, that is, pairs $(M, \alpha)$ where
  $M(\alpha)$ is not hyperbolic, of the types listed in
  Table~\ref{table:exceptsum} and distributed as in
  Figure~\ref{fig:numexcep}.
\end{theorem}
The list of these exceptional $(M, \alpha)$ together with the precise
topology of each $M(\alpha)$ is available at \cite{ExcepPaperData}.  Here,
in addition to describing the proof of Theorem~\ref{thm:except} in
Section~\ref{sec:proof}, I will give summaries of this data as it
relates to known results and open questions about Dehn filling in
Sections \ref{sec:exconj} and \ref{sec:new}.

\subsection{Prior work}

In the 1990s, Hodgson and Weeks studied the exceptional Dehn fillings
on the 286 manifolds in $\cC_5$; this work was never published but is
referred to extensively in \cite{Gordon1998} and provided many key
examples in the subject.  The series of papers
\cite{MartelliPetronio2006, MartelliPetronioRoukema2014, Martelli2018}
classified all exceptional fillings on an important series of chain
links with as many as 7 components; as noted in \cite[\S
3]{Martelli2018}, this determines the exceptional fillings on more
than 95\% of the 4{,}587 manifolds in $\cC_7$.  John Berge (personal
communication) independently did a search for exceptional fillings on
$\cC_9$ using a new version of his program Heegaard \cite{Heegaard3},
and found more than $99.3\%$ of the exceptional fillings included in
Theorem~\ref{thm:except}.

\subsection{Acknowledgements}

I thank Ken Baker, Bruno Martelli, and Jake Rasmussen for helpful
discussions on the topic of this paper and also thank John Berge for
sharing his unpublished data on this topic with me. This work was done
at the University of Illinois, the University of Melbourne, and IAS,
and was funded in part by the Simons Foundation and the US National
Science Foundation, the latter under grants DMS-1510204, DMS-1811156,
and the GEAR Network (DMS-1107452).  I also thank the referee for
their helpful comments on this paper, especially their observation in
Section~\ref{sec:slopelength}.

\section{Background and conventions}
\label{sec:conven}

\begin{figure}
\begin{center}
\begin{tikzpicture}[nmdstd,
  rounded corners=6,
  line width=1,
  every text node part/.style={align=center},
  boxed/.style = {rectangle, draw=black},
  outside sep/.style = {line width=1.2, dashed}
  ]
  \draw (0, 0) rectangle (3, 3);
  \draw (0, 3) rectangle (3, 6);
  \draw (7, 3) rectangle (10, 6);
  \draw[fill=nmddark!10]
  (3, 0) -- (10, 0) -- (10, 3) -- (7, 3) --
  (7, 6) -- (3, 6) --  cycle;
  \draw[fill=nmddark!20] (3, 0) -- (10, 0) -- (10, 3) --
  (7, 3) -- (7, 4) -- (3, 4) -- cycle;
  
  \begin{scope}[rounded corners=5, shift={(6, 1)}]
    \node[boxed] at (0.5, 0) {$S^3$};
    \node at (-1, 0) {lens \\ spaces};
    \node at (-0.5, 0.9) {finite $\pi_1$};
    \draw (-2, -0.5) rectangle (1, 0.5); 
    \draw (-2.5, -0.75) rectangle (1.5, 1.25);
  \end{scope}

  \node[boxed] at (9, 1) {$S^2  \times S^1$};
  \node[boxed] at (1.5, 1) {$\RP^3 \connsum \RP^3$};
  \node[boxed] at (1.5, 4.5) {$T^3 \connsum L(3, 1)$};
  \node[boxed] at (5.5, 4.5) {$\Sol$ manifolds};
  \node at (5, 3) {Seifert fibered \\ spaces};
  \node at (8.5, 4.5) {hyperbolic \\ piece in JSJ};
  \node at (5, 5.5) {graph manifolds};

  \node[above=0.1] at (1.5, 6) {conn.~sum\vphantom{prime}};
  \node[below=0.1] at (1.5, 0) {conn.~sum\vphantom{prime}};
  \draw[outside sep] (3, -0.6) -- (3, 6.6); 
  \node[above=0.1] at (6.5, 6) {prime};
  \node[below=0.1] at (6.5, 0) {prime};
  \node[left=0.3, rotate=90, anchor=center] at (0, 4.5) {toroidal};
  \node[left=0.3, rotate=90, anchor=center] at (0, 1.5) {atoroidal};
  \draw[outside sep] (-0.6, 3) -- (10.6, 3); 
  \node[right=0.3, rotate=-90, anchor=center] at (10, 4.5) {toroidal};
  \node[right=0.3, rotate=-90, anchor=center] at (10, 1.5) {atoroidal};
\end{tikzpicture}
\end{center}
\caption{The different types of nonhyperbolic 3-manifolds,
  and hence different types of exceptional Dehn fillings. See
  Section~\ref{sec:conven} for definitions.}
\label{fig:excepvenn}
\end{figure}
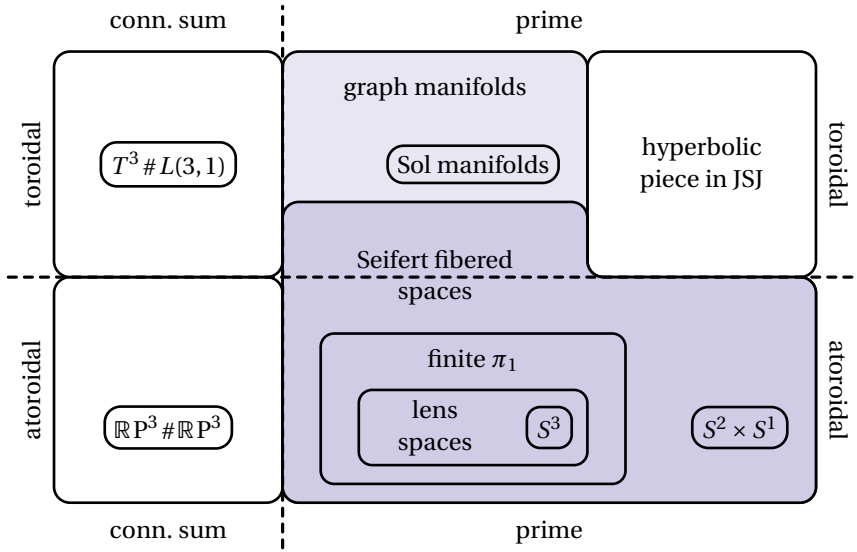

I first review the different types of nonhyperbolic \3-manifolds to
establish my conventions on the kinds of exceptional Dehn fillings one
can study.  Sources vary slightly on the latter point, and here I use
a relatively fine-grained division, which is illustrated in
Figure~\ref{fig:excepvenn}.  The summary for experts, who may safely
skip this section, is that here atoroidal means geometrically
atoroidal, the manifold $S^2 \times S^1$ is neither reducible nor an
honorary lens space, and the term Seifert fibered will not include
$\RP^3 \connsum \RP^3$.  I will assume familiarity with basic
3-manifold topology, the Geometrization Theorem, and the resulting
general structure of \3-manifolds, see e.g.~\cite{Hatcher3Manifolds,
  Scott1983, Bonahon2002} for details.  Throughout, all \3-manifolds
will be compact, orientable, and be either closed or have boundary
that is a union of tori; the symbol $M$ will always refer to such a
manifold.

Our first two kinds of nonhyperbolic \3-manifolds are those containing
certain spheres and tori.  An embedded \2-sphere in $M$ is
\emph{essential} if it does not bound a \3-ball.  If there are no
essential spheres then $M$ is \emph{irreducible}, and this includes
all hyperbolic $M$. Those $M$ containing \emph{separating} essential
spheres are called \emph{connected sums}; here, I avoid the more
common term reducible for this as for some authors reducible is the
complement of irreducible and so includes $S^2 \times S^1$ whose only
essential sphere is nonseparating.  When $M$ is not a connected sum it
is \emph{prime}. An $M$ is \emph{toroidal} when it contains an
embedded essential torus $T$, that is, one where $\pi_1 T \to \pi_1 M$
is injective and $T$ is not isotopic to a component of $\partial M$;
this is sometimes called \emph{geometrically toroidal}.  When $M$ is
not toroidal it is \emph{atoroidal}.  All hyperbolic $M$ are atoroidal.

When $M$ has a foliation by circles it is \emph{Seifert fibered} and
called a \emph{Seifert fibered space}; I will shortly revise this
definition to exclude a particularly unusual such manifold. The
Seifert fibered manifolds are exactly those admitting these six of the
eight possible geometries: $S^3$, $\E^3$, $\H^2 \times \R$,
$S^2 \times \R$, $\Nil$, and $\PSLRtilde$ and in particular are not
hyperbolic. Those with spherical geometry are precisely the
\3-manifolds with finite $\pi_1$, including $S^3$ itself and the
lens spaces $L(p, q)$ which are all quotients of $S^3$ by a cyclic
group. There are only two $M$ with an $S^2 \times \R$ geometry, namely
$S^2 \times S^1$ and $\RP^3 \connsum \RP^3$, which are both rather
special. First, $S^2 \times S^1$ is the only closed \3-manifold with Heegaard
genus one that is not a lens space, and also the only one whose
fundamental group is infinite cyclic. Second, $\RP^3 \connsum \RP^3$
is the unique Seifert fibered manifold that is a connected sum.  I
henceforth adopt the nonstandard convention that
$\RP^3 \connsum \RP^3$ is not Seifert fibered; this way, all Seifert
fibered spaces are prime.

Any irreducible $M$ has a collection of disjoint essential tori that cut
it up into pieces that are either Seifert fibered or hyperbolic.  The
minimal such collection is unique up to isotopy and gives the
\emph{JSJ decomposition} of $M$.  A \emph{graph manifold} is one where
all the pieces in the JSJ decomposition are Seifert fibered. So
Seifert fibered manifolds are graph manifolds as are those that admit
the $\Sol$ geometry; the latter are virtually torus bundles over the
circle with Anosov monodromy.  An irreducible $M$ that is neither
hyperbolic nor a graph manifold has a nontrivial JSJ decomposition
where at least one piece is hyperbolic; such $M$ have a
\emph{hyperbolic piece}.

Figure~\ref{fig:excepvenn} summarizes all the different types of
nonhyperbolic \3-manifolds.  Of course, many manifolds satisfy several
of these conditions, and in certain tables I will want each
nonhyperbolic manifold to have a single type.  In such cases, the type
used will be the most restricted possible in
Figure~\ref{fig:excepvenn}; for example, the type of $L(3, 1)$ will be
a lens space, even though it also has finite $\pi_1$, is Seifert
fibered, and is a graph manifold.  (This always makes sense
because I set up the various definitions to minimize overlaps that are
not containments.)  This more restrictive convention is used in
Tables~\ref{table:exceptsum} and Tables~\ref{table:exlens} only;
Table~\ref{table:distance} is correct with either convention.

\section{Evidence for standard conjectures}
\label{sec:exconj}

A great deal has been proven about the possibilities for exceptional
Dehn fillings; with regards to the gaps in our knowledge, the fillings
of Theorem~\ref{thm:except} are consistent with the standard
conjectures as I now describe.

\subsection{Knots in the 3-sphere} I start with the 1{,}267 manifolds in
$\cC_9$ that are exteriors of knots in $S^3$, which collectively have
some 2{,}615 additional exceptional fillings.

\begin{enumerate}[labelindent=-0.2em, labelsep=0.2em, leftmargin=*]
\item \label{item:cable} There are no Dehn fillings that are connected
  sums, consistent with the Cabling Conjecture \cite[\S
  2.2]{Gordon2003}.

\item The Berge Conjecture \cite[\S 3.2]{Gordon2003} holds for the 178
  nontrivial lens space fillings.
  
\item \label{item:intsur}
  All 1{,}143 nontrivial Seifert fibered fillings are along
  integral slopes and have the form $S^2(q_1, q_2, q_3)$ or
  $\RP^2(q_1, q_2)$; compare \cite[\S 3.3]{Gordon2003}.  In particular,
  all fillings with finite fundamental group are integral.

\end{enumerate}
We now turn to considering all of the manifolds in $\cC_9$.

\subsection{Distances between exceptional slopes}

A key topological invariant of a pair of slopes $\alpha$ and $\beta$
on a torus is their geometric intersection number
$\Delta(\alpha, \beta)$.  When $\alpha$ and $\beta$ are exceptional
slopes for a particular $M$, then $\Delta(\alpha, \beta) \leq 8$ by
\cite{LackenbyMeyerhoff2013}.  Gordon conjectured there are only four
possible $M$ with exceptional slopes where $ \Delta \geq 5$
\cite[Conjecture~3.4]{Gordon1998}, and this holds for $\cC_9$.  Much
of the work on exceptional fillings has focused on understanding the
maximum possible $\Delta(\alpha, \beta)$ where $M(\alpha)$ and
$M(\beta)$ are particular types of exceptional fillings. I summarize
what is observed for $\cC_9$ and how it relates to the known upper
bounds on $\Delta(\alpha, \beta)$ in Table~\ref{table:distance}.  In
all cases, the maximum value of $\Delta(\alpha, \beta)$ for $\cC_9$ is
the same as that already found in the literature, compare with
\cite{Gordon1998, Gordon1999, Gordon2003, Gordon2012} and also page
971 and Section A.2 of \cite{MartelliPetronio2006}.  I think it very
likely that all possible maximum values of $\Delta(\alpha, \beta)$
have been observed at this point.
  
\subsection{Atoroidal Seifert fibered and finite $\pi_1$ fillings}

There are two cases of $M$ in $\cC_9$ with slopes $\alpha$ and $\beta$
with $\Delta(\alpha, \beta) = 4$ where $M(\alpha)$ is an atoroidal
Seifert fibered space with $\pi_1\big(M(\alpha)\big)$ infinite and
$\pi_1\big(M(\beta)\big)$ is finite and noncyclic. These are already
contained in \cite{MartelliPetronio2006}, but this aspect is not
highlighted there and so is worth describing here.  The first example
is $m007$ where $m007(-2, 1)$ is the Seifert fibered space
$S^2\big((2,1), (3,1), (9,-7)\big)$ and $m007(2, 1)$ is
$S^2\big((2,1), (3,2), (3,-1)\big)$ which has nonabelian fundamental
group of order 120; this example is $M3_2$ with slopes $-4$ and $0$ in
Table A.3 of \cite{MartelliPetronio2006}.  The second is $m034$ with
$m034(2, 1) = S^2\big((2,1), (3,1), (11,-9)\big)$ and
$m034(-2, 1) = S^2\big((2, 1), (3, 2), (5, -3)\big)$ where the latter
has nonabelian fundamental group of order 2{,}040; it is the example
described in Table A.8 \cite{MartelliPetronio2006}, with $r/s = 2$ and
slopes $-4$ and $0$.  Here, my conventions for describing Seifert
fibered spaces follow Regina \cite{Regina}.

\subsection{Many exceptional fillings}
For a 1-cusped manifold $M$, let $e(M)$ denote the number of
exceptional fillings.  The distribution of $e(M)$ is shown in
Figure~\ref{fig:numexcep}.  There are only 11 manifolds in $\cC_9$
where $e(M) \geq 7$, namely $m003$, $m004$, $m006$, $m007$, $m009$,
$m016$, $m017$, $m023$, $m035$, $m038$, and $m039$.  According to
Gordon \cite[pages 136--7]{Gordon1998}, these 11 were first noticed
by Hodgson when he examined the 286 manifolds in $\cC_5$.  Gordon
writes there that ``In view of this data it is tempting to believe
that these eleven manifolds are the only ones with $e(M) \geq 7$'',
and it was later shown \cite {LackenbyMeyerhoff2013} that one always
has $e(M) \leq 10$.  These 11 are also the only manifolds with
$e(M) \geq 7$ among all Dehn fillings on the magic manifold
\cite{MartelliPetronio2006}.  In light of the additional data here,
it is safe to promote this temptation to a conjecture.

\subsection{Connected sums}

The connected sums in this census are all built of quite simple
pieces.  Specifically, the summands all have finite $\pi_1$ or are
$S^2 \times S^1$; there are only two summands in all but three cases:
the filling $o9_{39343}(1, 0)$ is
$\RP^3 \connsum \RP^3 \connsum \RP^3$ and both $o9_{41447}(1, 0)$ and
$o9_{43255}(1, 0)$ are the manifold
$L(3, 1) \connsum \RP^3 \connsum \RP^3$.  While there are infinite
families with two connected sum fillings \cite{Eudave-MunozWu1999}, it
is an open question whether there is a manifold with three such
fillings, see \cite[\S 4]{HoffmanMatignon2003}.  In $\cC_9$ there are
only 14 manifolds with two distinct Dehn fillings that are connected
sums, and none with more than two.  Another question from \cite[\S
4]{HoffmanMatignon2003} is when there are two such fillings, must both
have at least one summand that is $\RP^3 = L(2, 1)$, $L(3, 1)$, or
$L(4, 1)$?  The answer is yes for the 14 such manifolds in $\cC_9$.

\begin{table}
  \begin{center}
    \small
    \newcommand{\cgrey}{\cellcolor[gray]{0.8}}
    \def\arraystretch{1.7}
    \begin{tabular}{l|ccccccc}
      \toprule
       & $S^3$ & lens & finite & \#sum & $S^2 \times S^1$ & SFS(A) & tor \\

      \midrule

      $S^3$ &
      $-\infty$ &
      1 &
      1 \cgrey &
      $-\infty$ \cgrey &
      $-\infty$ &
      1 \cgrey &
      2 \\

      lens space &
      \cite{CullerGordonLueckeShalen1987} &
      1 &
      2 &
      1 &
      1 &
      2 \cgrey &
      3 \cgrey \\

      finite $\pi_1$ &
      \S \ref{sec:exconj}(\ref{item:intsur}) &
      \cite{BoyerZhang1996} &
      3 &
      1 &
      1 &
      4 \cgrey &
      5 \cgrey \\

      conn.~sum &
      \S \ref{sec:exconj}(\ref{item:cable}) &
      \cite{BoyerZhang1998} &
      \cite{BoyerGordonZhang2009} &
      1 &
      $-\infty$ \cgrey &
      2 \cgrey & 3 \\

      $S^2 \times S^1$ &
      \cite{Gabai1987c} &
      \cite{CullerGordonLueckeShalen1987} &
      \cite{BoyerGordonZhang2009} &
      \cite{GordonLuecke1996}&
      $-\infty$ &
      1 \cgrey &
      2  \\

      Seifert (ator.) &
      \S \ref{sec:exconj}(\ref{item:intsur}) &
      &
      &
      \cite{BoyerCullerShalenZhang2008} &
      &
      6 \cgrey &
      7 \cgrey \\ 

      toroidal &
      \cite{GordonLuecke1995} &
      \cite{Lee2011} &
      &
      \cite{Oh1997, Wu1998} &
      \cite{Lee2005} &
      \cite{BoyerGordonZhang2012}
      & 8 \\
      
      \bottomrule
    \end{tabular}
    \caption{This table lists the maximum value of
      $\Delta(\alpha, \beta)$ where $M \in \cC_9$ and the fillings
      $M(\alpha)$ and $M(\beta)$ have the indicated types; here it is
      implicit that $\alpha \neq \beta$ and $-\infty$ is used when no
      such pair exists. As these values form a symmetric matrix, only
      the entries on or above the diagonal are given.  Below the
      diagonal is a reference for the strongest proven upper bound on
      $\Delta(\alpha, \beta)$ for that pair of types, with references
      for the diagonal being, in order, \cite{GordonLuecke1989,
        CullerGordonLueckeShalen1987, BoyerZhang2001,
        GordonLuecke1996, Gabai1987b, LackenbyMeyerhoff2013,
        Gordon1998TAMS}.  Cases where the value for $\cC_9$ is
      \emph{strictly smaller} than the best proven upper bound are
      indicated by the shaded boxes.  Thus in the unshaded cases the
      maximum possible value of $\Delta(\alpha, \beta)$ for \emph{all}
      \3-manifolds has been established exactly by the indicated
      reference.  }\label{table:distance}
  \end{center}
       
\end{table}

\section{New observations}
\label{sec:new}

Here are some interesting patterns that I couldn't find in the
existing literature.  I encourage you to download the complete
data at \cite{ExcepPaperData} and find others that I have missed.

\subsection{Finite nonabelian fillings} The maximum number of fillings
on $M$ in $\cC_9$ where the fundamental group is finite and nonabelian
is three.  There are only four such $M$, namely $m011$, $s757$,
$v2702$, and $v2797$.  I conjecture that these are the
only four manifolds with this property.
  
\subsection{Toroidal fillings} The maximum number of toroidal fillings
on $M$ in $\cC_9$ is 4, and there only 27 such $M$, namely $s772$,
$s778$, $s911$, $v2640$, $t08282$, $t11538$, $t12033$, $t12035$,
$t12036$, $t12041$, $t12043$, $t12045$, $t12050$, $t12548$, $t12648$,
$o9_{35259}$, $o9_{36732}$, $o9_{37030}$, $o9_{38039}$, $o9_{39094}$,
$o9_{40054}$, $o9_{41000}$, $o9_{41004}$, $o9_{41006}$, $o9_{41007}$,
$o9_{41008}$, $o9_{43799}$. Are there are infinitely many such
examples?  Perhaps we should expect there to be since the previous
list includes manifolds with 6, 7, 8, and 9 ideal tetrahedra.  None of
the examples with 4 toroidal Dehn fillings is the exterior of a knot
in $S^3$, consistent with a conjecture of \cite[Page
60]{Eudave-Munoz1997}.

\begin{table}
  \begin{center}
    \small
    \begin{tabular}{lcc}
      \toprule
      \multicolumn{1}{c}{Type of $M(\alpha)$} &
      $\min\big(\ell(\alpha)\big)$ & $\max\big(\ell(\alpha)\big)$ \\
      \midrule
      $S^3$              & 1.000 & 3.323  \\
      $S^2 \times S^1$   & 1.288 & 3.328  \\ 
      connected sum      & 1.398 & 3.707  \\ 
      lens space         & 1.189 & 3.928  \\ 
      Sol torus bundle   & 2.288 & 4.185 \\ 
      finite $\pi_1$     & 1.520 & 4.443  \\ 
      Seifert (toroidal) & 1.906 & 4.583  \\ 
      Seifert (atoroidal)& 1.935 & 4.841  \\
      graph manifold     & 2.178 & 5.318  \\ 
      hyperbolic piece   & 3.520 & 6.000 \\
      \bottomrule
    \end{tabular}
  \end{center}
  \caption{This table gives minimum and maximum lengths of each type
    of exceptional slope in Theorem~\ref{thm:except} as measured in
    the torus bounding a maximal cusp. Lengths have been rounded to
    three decimal places, and here type refers to the most restricted
    category possible from Figure~\ref{fig:excepvenn}, which is why
    the $\min\big(\ell(\alpha)\big)$ for lens spaces is bigger than
    that for $S^3$. Compare with \cite[Table
    1]{HoffmanPurcell2017}.  } \label{table:exlens}
\end{table}

\begin{table}
  \begin{center}
    \small
    \begin{tabular}{rrrrc}
      \toprule
      \multicolumn{1}{c}{$M$} &    \multicolumn{1}{c}{$\alpha$} &
      \multicolumn{1}{c}{$M(\alpha)$} & $\abs{\pi_1\big(M(\alpha)\big)}$
      & \multicolumn{1}{c}{$\ell(\alpha)$} \\
      \midrule
         $t05002$ &  $(-1, 1)$ &   $S^2\big((2,1), (3,2), (3,-1)\big)$ & 120 &  4.004139 \\
      $o9_{28194}$ &  $(-1, 1)$ &   $S^2\big((2,1), (3,2), (4,-1)\big)$ & 528 &    4.017192 \\
      $o9_{35417}$ &   $(1, 1)$ &    $S^2\big((2,1), (2,1), (8,3)\big)$ & 352  &   4.021047 \\
      $o9_{35418}$ &   $(1, 1)$ &  $S^2\big((2,1), (2,1), (11,-3)\big)$ & 352  &   4.021047 \\
          $v3479$ &  $(1, 1)$ &   $S^2\big((2,1), (3,2), (5,-1)\big)$ &  3{,}480 &   4.028619 \\
      $o9_{36221}$ &   $(1, 1)$ &   $S^2\big((2,1), (2,1), (10,3)\big)$ & 520  &   4.033399 \\
      $o9_{36224}$ &  $(-1, 1)$ &  $S^2\big((2,1), (2,1), (13,-3)\big)$ & 520  &   4.033399 \\
         $v2420$ &  $(-1, 1)$ &   $S^2\big((2,1), (3,2), (4,-1)\big)$ & 528  &   4.060890 \\
           $m342$ &   $(1, 1)$ &   $S^2\big((2,1), (3,2), (3,-1)\big)$ & 120  &   4.067597 \\
          $m011$ &   $(2, 1)$ &   $S^2\big((2,1), (3,2), (5,-3)\big)$ & 2{,}040 &    4.085768 \\
      $o9_{41134}$ &  $(-1, 1)$ &   $S^2\big((2,1), (3,2), (4,-1)\big)$ & 528 &    4.184451 \\
     $o9_{12592}$ &  $(-1, 1)$ &   $S^2\big((2,1), (3,2), (5,-3)\big)$ &  2{,}040 &   4.195283 \\
          $s954$ &   $(1, 1)$ &   $S^2\big((2,1), (3,2), (4,-1)\big)$ &  528 &   4.207000 \\
          $s546$ &  $(-1, 1)$ &   $S^2\big((2,1), (3,2), (5,-3)\big)$ &  2{,}040 &   4.442966 \\
      \bottomrule
    \end{tabular}

    \vspace{-1em}
  \end{center}
    \caption{All pairs $(M, \alpha)$ with $M \in \cC_9$ and
      $\pi_1\big(M(\alpha)\big)$ finite where $\ell(\alpha) > 4$.}
    \label{table:longslopes}

\end{table}

\subsection{Lengths of exceptional slopes}
\label{sec:slopelength}
Hoffman and Purcell \cite{HoffmanPurcell2017} studied the length of
exceptional slopes $\alpha$ in the horotorus cutting off a maximal
cusp for $M$.  By the $6$-Theorem, the length $\ell(\alpha)$ of such
$\alpha$ is at most $6$. Table~\ref{table:exlens} details the longest
exceptional slopes of each type observed in $\cC_9$; compare with
Table 1 of \cite{HoffmanPurcell2017}. The new feature is slopes of
length more than 4 yielding manifolds with finite fundamental group;
these are listed in Table~\ref{table:longslopes}.  Can some of these
be made into an infinite family of finite exceptional slopes with
$\ell(\alpha) \to 5$ analogous to Proposition~4.2 of
\cite{HoffmanPurcell2017}?

The referee kindly pointed out that one can use a covering trick to
create even longer slopes for some types starting with the examples in
Table~\ref{table:longslopes}.  Specifically, the extreme example for
lens spaces comes from $M = o9_{18855}$ and the slope
$\alpha = (1, 1)$ where $M(\alpha) = L(39,16)$ and
$\ell(\alpha) \approx 3.92794$.  The core curve of the Dehn filling
turns out to generate $H_1(M(\alpha); \Z) \cong \Z/39\Z$, and
consequently one can take a $39$-fold cyclic cover of $M$ to get the
exterior of a knot in $S^3$ whose meridian $\mu$ also has
$\ell(\mu) \approx 3.92794$. As discussed in \cite{HoffmanPurcell2017},
it is conjectured that for an $S^3$ filling one always has
$\ell(\mu) \leq 4$, and there are several families of such where
$\ell(\mu) \to 4$ from below.

One can apply the same trick to $s546(-1, 1)$ from
Table~\ref{table:longslopes} to produce a hyperbolic knot in the
Poincar\'e homology sphere where the meridian has length about
$4.442966$; specifically, take the 17-fold cyclic cover of $s546$
corresponding to the kernel of the map
$\pi_1\big(s546(-1, 1)\big) \to H_1\big(s546(-1, 1); \Z\big) \cong
\Z/17\Z$.

\subsection{A cabling conjecture for $S^2 \times S^1$} As per
Table~\ref{table:distance}, there are no known hyperbolic knot
exteriors in $S^2 \times S^1$ with a Dehn filling that is a connected
sum. Thus, as in Section \ref{sec:exconj}(\ref{item:cable}) for the
case of $S^3$, I conjecture that none exist, i.e.~that
no hyperbolic knot in $S^2 \times S^1$ has a Dehn surgery yielding a
connected sum.

\section{Outline of the proof of Theorem~\ref{thm:except}}
\label{sec:proof}
  
I turn now to the proof of Theorem~\ref{thm:except}.  Initially, I
found a candidate $\cE$ for the list of all exceptional fillings as a
byproduct of another project.  However, the proof of the correctness
of $\cE$ follows the approach of \cite{MartelliPetronioRoukema2014}.
The list $\cE$, related data, and the code used in the proof can all
be obtained from \cite{ExcepPaperData}; to run the code, which
requires using several software packages together in consort, the
Docker image \cite{KitchenSink} may be helpful.

\begin{proof}[Proof of Theorem~\ref{thm:except}]
  The set $\cE$ consists of 205{,}822 pairs $(M, \alpha)$ where
  $M \in \cC_9$ and $\alpha \in \slopes(M)$. There are two things to
  show: that every $M(\alpha)$ is not hyperbolic of the type
  claimed in Table~\ref{table:exceptsum} and that all other fillings
  on $M \in \cC_9$ are hyperbolic.

  For the latter task, for each $M \in \cC_9$ I found an embedded cusp
  neighborhood so that I could measure the lengths of slopes in its
  horotorus boundary; this was done rigorously in SnapPy \cite{SnapPy}
  running inside \cite{SageMath} using the approach of
  \cite{HIKMOT2016} and \cite[\S 3.6]{DunfieldHoffmanLicata2015}.  By
  the 6-Theorem of \cite{Agol2000, Lackenby2000}, it suffices to
  examine all slopes $\beta \in \slopes(M)$ where
  $\ell(\beta) \leq 6$.  For the cusp neighborhoods I used, overall
  there were some 355{,}128 such slopes.  For the 149{,}306 pairs
  $(M, \beta)$ that were not in $\cE$, I checked that $M(\beta)$ was
  hyperbolic using the method of \cite{HIKMOT2016} as reimplemented in
  SnapPy.  As in the proof of Theorem 5.2 of \cite{HIKMOT2016}, it was
  sometimes necessarily to search around for a triangulation that
  could be used to certify the existence of a hyperbolic structure.
  This completes the proof that filling along any slope not in $\cE$
  yields a hyperbolic manifold.
  
  In the other direction, to show that each $M(\alpha)$ in $\cE$
  is not hyperbolic, I primarily used Regina \cite{Regina},
  specifically its combinatorial recognition methods \cite[\S
  4]{Burton2013}.  These work when the input triangulation has the
  very particular form associated to a standard triangulation of a
  Seifert fibered space or graph manifold. Of course, there are many
  triangulations of such manifolds which do not have this structure,
  so I generated many different 1-vertex triangulations of each
  $M(\alpha)$ and fed them into Regina until it succeeded in
  recognizing the topology.  This worked for all but 2{,}890 of the
  $M(\alpha)$.  Of those remaining, in 680 cases the Recognizer
  program of \cite{Matveev1998, Spine} showed that they were graph
  manifolds.  (Currently, Regina can only identify graph manifolds
  where the graph in question is either a segment with two or three
  vertices or a loop with one vertex. These 680 all have slightly more
  complicated graphs, for example a loop with either two or three
  vertices.)  For each of the remaining 2{,}210 manifolds, Regina
  found at least one essential normal torus.  Cutting along a suitable
  collection of such essential tori gave pieces that always included a
  cusped hyperbolic 3-manifold with an ideal triangulation with at
  most 6 ideal tetrahedra; in particular, each of these 2{,}210
  manifolds is nonhyperbolic with a non-trivial JSJ decomposition with
  a hyperbolic piece.  Thus every $M(\alpha)$ in $\cE$ is not
  hyperbolic.  This completes the proof that $\cE$ is precisely the
  list of exceptional fillings on the manifolds in $\cC_9$.

  To prove the correctness of Table~\ref{table:exceptsum}, the hard
  part is ensuring that nothing listed as a (proper) graph manifold is
  actually Seifert fibered.  Everything else can be read off from the
  Seifert/graph descriptions found in the previous step, though I
  double-checked much of it in other ways.  For example, I used Magma
  \cite{Magma} to give an independent check that the 59{,}200
  spherical manifolds had the claimed type of fundamental group (this
  could also be done with GAP \cite{GAP}).  I also had Regina compute
  directly which manifolds are toroidal using normal surface
  techniques and this matched what follows from the Seifert/graph
  descriptions.  As mentioned, Regina identifies structure in the
  given triangulation, which might well be a graph manifold structure
  that can be simplified after the fact. For example, it will
  sometimes return graph manifolds where one of the nodes is a solid
  torus.  In such instances, additional triangulations were examined
  until a more concise description was found.  To certify a graph
  description as minimal, I just checked that all Seifert pieces have
  incompressible boundary (i.e.~no solid tori) and that no two Seifert
  pieces are glued together so that the fibers match up; here, I took
  care to consider the possibility of switching the Seifert fibration
  for the exceptional piece which is both
  $D^2\big( (2, 1), (2, 1)\big)$ and the twisted circle bundle over
  the M\"obius band.
\end{proof}

{\RaggedRight
\small    
\bibliographystyle{nmd/math} 
\bibliography{conjecture}
}
\end{document}